# A MARKOV CHAIN MODEL OF A POLLING SYSTEM WITH PARAMETER REGENERATION

By Iain MacPhee, Mikhail Menshikov,[1] Dimitri Petritis[2] and Serguei Popov [3]

*University of Durham, University of Durham, Université de Rennes 1 and Universidade de São Paulo*

We study a model of a polling system, that is, a collection of $d$ queues with a single server that switches from queue to queue. The service time distribution and arrival rates change randomly every time a queue is emptied. This model is mapped to a mathematically equivalent model of a random walk with random choice of transition probabilities, a model which is of independent interest. All our results are obtained using methods from the constructive theory of Markov chains. We determine conditions for the existence of polynomial moments of hitting times for the random walk. An unusual phenomenon of thickness of the region of null recurrence for both the random walk and the queueing model is also proved.

**1. Introduction.** The results presented in this paper, depending on the affinities of the reader, can be looked at from two different viewpoints: a model from queueing theory or a model from the theory of Markov chains. These two models are mathematically equivalent but have their own individual interest. In this paper we start by describing the queueing model and stating the results for it. Subsequently, the proofs are obtained, after the queueing model has been bijectively mapped to a random walk model, by using methods from the constructive theory of Markov chains to determine the asymptotic behavior of this random walk. A phenomenon, exceptional in

Received July 2006; revised December 2006.
[1]Supported in part by FAPESP (20004/13610-2).
[2]Supported in part by the "Réseau Mathématique France–Brésil" and by the European Science Foundation.
[3]Supported in part by CNPq (302981/02–0), FAPESP (04/03056–8) and by the "Rede Matemática Brasil–França."

*AMS 2000 subject classifications.* Primary 60K25, 60J10; secondary 60G42, 90B22.
*Key words and phrases.* Polling system, parameter regeneration, stability, time-inhomogeneous Markov chains, recurrence, Lyapunov functions.

the theory of stochastic processes, namely, a thick region of null recurrence, is also proved for our model.

Our main objective is to investigate stability conditions for a polling model with *parameter regeneration*. We feel that parameter regeneration is quite a natural hypothesis for realistic queueing systems: it models changes in behavior (of customers or servers) due to factors that are *external* to the studied system. In our model there are $d \geq 2$ stations where queues of jobs can form and a single server which switches from station to station. The server only processes jobs at its current station. For clearer exposition we treat in this paper the case where only two stations are *open*, that is, receiving jobs, at any time. We study the model where all stations are open in a separate paper [9]. The reason for splitting this work into two separate papers, depending on whether two or more stations receive jobs, is that in the former case stability conditions are expressible in terms of *explicitly* determined factors, while in the latter, they are expressible in terms of factors whose existence is shown but whose values are inaccessible to direct computation. These factors are reminiscent of the Lyapunov exponents for products of random matrices whose values can be approximated only by stochastic simulation. Moreover, while the obtained results in the two papers can be phrased in similar ways, the techniques used to obtain them differ at key points.

When the server completes all the jobs at its current station $k$, say, it starts to move to the other open station $i$, station $k$ closes and independently of the process history, station $j \neq i$ opens with probability $P_{ij}$, where $\sum_{j \neq i} P_{ij} = 1$ for each station $i$. The sequence of stations visited by the server thus forms a Markov chain $\{I_r : r = 1, 2, \ldots\}$ on state space $\mathcal{S}_d = \{1, 2, \ldots, d\}$ with transition matrix $P$ which we will assume is irreducible but may be periodic (this is necessary for the case $d = 2$ and for cases where one station has priority over the others). We call the time period between the server's arrival at and departure from a station an *epoch*. The server takes some time with distribution $H_{ij}$ which is independent of the process history to switch between stations and we will assume such times have bounded first moments (for some results a stronger assumption will be needed).

During each switching time and service epoch, jobs arrive at the open stations in independent Poisson streams with rates $\lambda_1$ at the server's station $i$ and $\lambda_2$ at the other open station $j$. The jobs at station $i$ have independent service times with distribution function $G$, mean $\mu^{-1}$, where $G$ is selected from the class of distributions $\Gamma$ with uniformly bounded second moments. This service time distribution and the arrival rates $(\lambda_1, \lambda_2)$ are regenerated at the start of each switching time independently of the process history using a measure $\nu_{ij}$ which depends upon the ordered pair $(i, j)$ of open stations. It will follow from our Condition E, given later in the paper, that there exist two positive constants $m_0$ and $M_0$, with $m_0 < M_0$ such that $m_0 < \lambda_i < M_0$ for $i = 1$, 2. Each job at queue $i$ leaves the system after processing by the server and when the server completes all jobs at station $i$, including any that arrive since the server switched to $i$ (exhaustive service), the epoch ends. Succeeding epochs are defined in the same way except that if the system ever has no jobs, then the server waits for some to arrive—we do not go into any details as we are interested in establishing stability criteria, not studying equilibrium behavior. Our main interest is in the queue lengths, but their behavior is controlled by the process of regenerated parameters, so we consider this in detail.

*Parameter regeneration.* The parameters in each epoch have the form $(i, j, G, \lambda_1, \lambda_2) \in \mathcal{S}_d^2 \times \Gamma \times (m_0, M_0)^2 \equiv \Omega_e$. We impose below conditions on service times and arrival rates that ensure that the overall rate of events for the random walk is uniformly bounded and all epochs are finite with probability 1.

As stated above, we assume that the Markov chain $\{I_r\}$ supplying the sequence of open stations is irreducible and starts from a state $i_1 \in \mathcal{S}_d$. Denote the space of parameter sequences by $\Omega_S = \Omega_e^{\mathbb{N}} = \{(i_r, i_{r+1}, G_r, \boldsymbol{\lambda}_r) : r = 1, 2, \ldots\}$, the collection of all sequences of ordered pairs of open stations together with service time distribution and arrival rates.

Probability is associated with events in $\Omega_S$ as follows. At each switching time $r \geq 1$, the station where the server was closes, the server starts switching to station $I_r = i$ and at the same time, with chance $P_{ij}$, $j = 1$, ..., $d$, station $I_{r+1} = j \neq i$ opens. The arrival rates $\boldsymbol{\lambda}_r = (\lambda_1(r), \lambda_2(r))$ and the service time distribution $G_r$ at queue $i$ (the *random parameters*) are chosen using measure $\nu_{ij}$ independently of the process history and are constant throughout the epoch. As stated above, the server takes a time with distribution $H_{ij}$ to complete its switch from $i$ to $j$.

We will use $\omega \in \Omega_S$ to denote a sequence of regenerated parameters and $\omega_r = (i_r, i_{r+1}, G_r, \boldsymbol{\lambda}_r)$ to denote the parameter set in epoch $r$. Where it will not cause confusion we will drop the epoch number dependence from parameter sets and write $\omega_r = (i, j, G, \boldsymbol{\lambda})$.

*Embedded Markov chain.* We find it convenient to establish our stability results for the queueing process in a discrete time setting, so we work with a Markov chain embedded in the continuous time process. Specifically, we consider the discrete time embedded process $(W, \xi) = \{W_n, \xi_n\}_{n \in \mathbb{N}}$ by observing the continuous time process just after the instants of service completions and the ends of switching times. Here $W_n \in \Omega_e$ identifies the parameters immediately after event $n$ and $\xi_n$ denotes the ordered pair of queue lengths so $\xi_n \in \mathbb{N}_0^2$ (where $\mathbb{N}_0 = \{0, 1, 2, 3, \ldots\}$) and specifically, $\xi_{n,1}$ is the queue length at the server's current station and $\xi_{n,2}$ the queue length at the other open station.

Let $\mathbf{e}_1$ and $\mathbf{e}_2$ denote the unit vectors of $\mathbb{R}^2$. Given $W_n = (k, i, G, \boldsymbol{\lambda})$, the possible increments of the queue lengths at event $n$ are of the form

$$s\mathbf{e}_1 + t\mathbf{e}_2, \qquad s, t \in \mathbb{N}_0 \text{ during a switching time}$$

and

$$s\mathbf{e}_1 + t\mathbf{e}_2, \qquad s \in \{-1\} \cup \mathbb{N}_0; t \in \mathbb{N}_0 \text{ during a service epoch}$$

(as the server's queue is always given by $x_1$). During switching times and during service epochs when $\xi_n = (x_1, x_2)$ with $x_1 > 1$, the parameters will not change, that is, $W_{n+1} = W_n$ but if, during a service epoch, we have $x_1 = 1$, then an increment $(-1, t)$ ends the current epoch. For the continuous time process, as stated earlier: station $k$ closes; the server starts the switch to station $i$; some station $j$, selected using the Markov chain $\{I_r\}$, opens; given the pair $(i, j)$ a service time distribution $G$ and arrival rates $\boldsymbol{\lambda}$ are selected using measure $\nu_{ij}$; the customers waiting at station $i$ become the initial queue at the server's new location and during the switching time customers arrive at the open stations $(i, j)$ with the selected rates. The possible transitions for our embedded process are to states

$$(i, j, \hat{G}, \hat{\boldsymbol{\lambda}}; x_2 + s, t) \quad \text{with intensity} \quad P_{ij} d\nu_{ij}(\hat{G}, \hat{\boldsymbol{\lambda}}) h_{ki}(s, t; \boldsymbol{\lambda}),$$

where $h_{ki}(s, t; \boldsymbol{\lambda})$ is the probability of increment $s\mathbf{e}_1 + t\mathbf{e}_2$, given the arrival rates $\boldsymbol{\lambda}$, at the transition where the server switches from station $k$ to $i$.

In order to study the stability of the queueing system, we will consider the trajectories of a particle moving according to drift vector fields $\mathbf{D}$ related to the parameter regeneration process. As such vector field or *fluid* models are quite standard in the queueing literature, we will leave details for later and briefly describe now only what we need to state our main results. The *one-step mean drifts* of $\xi$, conditioned on the parameters in the service epochs, form a vector field $\mathbf{D}^0 : \Omega \times \mathbb{N} \to \mathbb{R}^2$ defined by

$$(1) \qquad \mathbf{D}^0(\omega, r) = \mathbf{E}(\xi_{n+1} - \xi_n \mid \omega_r) = \left( \frac{\lambda_1}{\mu} - 1, \frac{\lambda_2}{\mu} \right),$$

where $\omega_r = (i, j, G, \lambda_1, \lambda_2)$. The process $\xi$ is homogeneous during each service epoch so the mean drift does not depend on queue lengths $\xi_n$ but only upon $\omega_r = (i_r, i_{r+1}, G_r, \boldsymbol{\lambda}_r)$, the parameters in epoch $r$. Our vector field model will not include switching times so we only define the mean drifts for service epochs.

Now we observe that if a particle starts from position $(x, 0)$ and moves across a quarter plane in direction $\mathbf{D}^0$, then it eventually arrives at point $(0, x\lambda_2/(\mu - \lambda_1))$. The multiplying ratio that appears here plays a key role in the phenomena we are interested in. For later use we observe that, given open stations $(i, j)$,

$$(2) \qquad Y_{ij} \equiv \frac{\lambda_2}{\mu - \lambda_1} = \tan \theta,$$

where $\theta$ is the minor angle between $\mathbf{D}$ and the vector $-\mathbf{e}_1$. We also introduce the notation

$$F_{ij}(x) = \nu_{ij}[(G, \lambda_1, \lambda_2) : Y_{ij} \leq x] \qquad \text{for } x > 0$$

for the cumulative distribution function of the slope $Y_{ij}$ for each pair of stations $i, j \in \mathcal{S}_d = \{1, 2, \ldots, d\}$ with $i \neq j$. For notational convenience, set $F_{ii}(x) = 0$, $x < 1$ and $F_{ii}(x) = 1$, $x \geq 1$, for all $i \in \mathcal{S}_d$.

*Regularity conditions.* In order to obtain our results, we have to impose some restriction on the regeneration of parameters and we state these now.

CONDITION E. For all pairs of distinct stations $i$, $j$, there exists $m_0 > 0$, $M_0 > m_0$ such that:

(i) $\nu_{ij}[(G, \boldsymbol{\lambda}) : \lambda_1 + m_0 < \mu < M_0] = 1$,
(ii) $F_{ij}(1/M_0) = 0$, $F_{ij}(M_0) = 1$.

Part (i) ensures that the service rate is uniformly bounded from above and below and that the server is never trapped at any station forever. Part (ii) limits the effect that any single epoch can have on the long term behavior of the process by preventing the mean drift vectors $\mathbf{D}(\omega, r)$ from approaching the axial directions too closely (so it is a type of uniform ellipticity condition).

By observing the vector field model only at parameter regeneration times we obtain a *multiplicative Markov process* of independent interest which we now define.

*Multiplicative model.* Let $(I, X) = \{(I_r, X_r) : r \geq 0\}$ denote a discrete time Markov process on the state space $\mathcal{S}_d \times \mathbb{R}_+$ (recall that $\mathcal{S}_d = \{1, 2, \ldots, d\}$) which jumps from state $(i, x) \mapsto (j, \alpha x)$ with transition measure $P_{ij} dF_{ij}(\alpha)$ for $\alpha > 0$. The matrix $P = [P_{ij}]$ is the transition matrix of an irreducible Markov chain on $\mathcal{S}_d$ and the $F_{ij}$ are a family of distribution functions with support $\mathbb{R}_+$ that satisfy Condition E(ii). Assuming that the process waits for time $x$ in state $(i, x)$ before making its next jump, we associate with $(I, X)$ an *accumulated time* process $\{T_r\}$ defined by

(3) $$T_0 = 0 \quad \text{and} \quad T_r = T_{r-1} + X_r = \sum_1^r X_n, \qquad r = 1, 2, \ldots.$$

The sequence of times $\{T_r\}$ is increasing and so the limit $T = \lim_{r \to \infty} T_r$ exists (in some cases it is infinite) and we call it the *total time* of the multiplicative model.

**2. The main results.** Let $\pi = (\pi_1, \pi_2, \ldots, \pi_d)$ denote the equilibrium distribution of the open stations Markov chain $\{I_r\}$. For each pair $i, j \in \mathcal{S}_d^2$, let $Y_{ij}$ be a random variable with distribution $F_{ij}$ and for any constant $s > 0$, set

$$(4) \qquad L_i = \sum_{j=1}^{d} P_{ij} \mathbf{E}(\log Y_{ij}) \quad \text{and} \quad M_{ij}(s) = P_{ij} \mathbf{E}(Y_{ij}^s).$$

The $d \times d$ matrix $M(s) = [M_{ij}(s)]$ is nonnegative and irreducible because $\{I_r\}$ is irreducible. Hence, $M(s)$ has a Perron–Frobenius eigenvalue $\eta(s) > 0$ such that $\eta(s) \geq |\eta_i(s)|$ for all other eigenvalues of $M(s)$. The size of $\eta(s)$ determines which moments of the total time $T$ of the multiplicative model exist, which in turn affects the long term behavior of the random walk $\xi$.

It is important in Theorem 3.1 below to note that

$$(5) \qquad \eta'(0) = \sum_{i=1}^{d} \pi_i L_i.$$

We briefly indicate why this is the case. We have $M'_{ij}(0) = P_{ij} \mathbf{E}(\log Y_{ij})$ and we note that the matrix $M(0)$ has $\eta(0) = 1$ with corresponding left eigenvector $\pi$ and right eigenvector $\mathbf{1}$. For any fixed $s > 0$, let $u(s)$ be the left eigenvector corresponding to eigenvalue $\eta(s)$ of $M(s)$. Scale $u(s)$ so that $(\pi - u(s)) \cdot \mathbf{1} = 0$ ($\cdot$ indicates the scalar product). Now, for $s$ small, write $u(s) = \pi + sx$, where $x \cdot \mathbf{1} = 0$. We note that $M_{ij}(s) = P_{ij}(1 + s\mathbf{E}(\log Y_{ij})) + O(s^2)$ and $\eta(s) = u(s) M(s) \mathbf{1} (u(s) \mathbf{1})^{-1} = u(s) M(s) \mathbf{1}$.

Now we expand $u(s) M(s) \mathbf{1}$ to get

$$\eta(s) = (\pi + sx)[P_{ij}(1 + s\mathbf{E}(\log Y_{ij}))]\mathbf{1} + O(s^2) = 1 + s \sum \pi_i L_i + O(s^2)$$

since $xP\mathbf{1} = x \cdot \mathbf{1} = 0$. The relation $\eta(s) = 1 + s \sum \pi_i L_i + O(s^2)$ immediately implies (5).

The multiplicative Markov chain's time component and our calculations for it are key elements of our main results which give stability conditions for the queue length component $\xi = \{\xi_n\}$ of the stochastic model $(W, \xi)$. Let $\varnothing$ denote the state $(0, 0)$ where both open queues have no customers and let $\tau = \min\{n > 0 : \xi_n = \varnothing\}$ be the time until the system first reaches this empty state. In both of the following theorems it is assumed that the measures $\nu_{ij}$ for selecting the parameters at the starts of epochs satisfy Condition E.

THEOREM 2.1.  (i) *If* $\sum_1^d \pi_i L_i < 0$, *then* $\mathbf{P}(\tau < \infty) = 1$.

(ii) *If* $\sum_1^d \pi_i L_i = 0$ *and switching times are all 0, then* $\mathbf{P}(\tau < \infty) = 1$.
(iii) *If* $\sum_1^d \pi_i L_i > 0$, *then* $\mathbf{P}(\tau = \infty) > 0$.

REMARK. The case where $\sum_1^d \pi_i L_i = 0$ and switching times are nonzero is difficult and we do not consider it here. The two station model in this case but with constant, not regenerated, parameters was considered by Menshikov and Zuyev in [11]. A detailed classification of a critical random walk on a 2-dimensional complex, again with a fixed jump distribution on each face, was given by MacPhee and Menshikov in [8]. In these papers, generalizations of the results of Lamperti on processes with asymptotically zero drift (see [6]) were used to establish how the behavior of the process during switching events affects which moments of hitting times to the empty state exist.

The major result in this paper is that we can make refined statements about the tail of the hitting time $\tau$ in the recurrent case (i) of the previous result by determining which moments of $\tau$ are finite and which are infinite. The existence of moment $s$ of $\tau$ requires moment conditions on the service and switching times, specifically that there is some constant $K_s$ such that $\int_0^\infty x^s \, dG(x) < K_s$ uniformly over $G \in \Gamma$ and, similarly, $\int_0^\infty x^s \, dH_{ij}(x) < K_s$ for all pairs of stations $i$, $j$.

THEOREM 2.2. *If $\sum_1^d \pi_i L_i < 0$, then:*

(i) *if $\eta(s) < 1$ and the sth moments of switching times and service times are uniformly bounded by some constant $K_s$, then $\mathbf{E}(\tau^s) < \infty$;*

(ii) *if $\eta(s) \geq 1$, then $\mathbf{E}(\tau^s) = \infty$.*

This result implies that the region of parameter space where the queueing process $\xi$ exhibits "null-recurrence" can be of considerable volume. It also implies that we have only polynomial ergodicity in the sense that the time of reaching the empty state has only finite order moments and the order smoothly changes when we pass from the null recurrent to the positive recurrent case. Note that as long as they have sufficient moments, the switching times play no role whatever in the result which is in accord with what is known about the constant parameter polling network in the recurrent case; see [5]. The phenomenon of a large null-recurrence region is present even in the case of zero switching times. We now illustrate it with a simple example before further commenting on it.

EXAMPLE. There are $d = 2$ nodes and the arrival rates are $\lambda_1 = \lambda_2 = 1$ in all epochs. When the server is at station 2 service times are always exponential with rate 2. At station 1 service times are again exponential, but there are two possible service rates, $5/4$ and $5$. Each time the server comes to station 1 it works at the smaller rate with probability $p \in [0, 1]$ and at the larger rate with probability $1 - p$. That is, according to our notation,

$$\nu_{12}[(\mu, \boldsymbol{\lambda}) : \mu = 2, \boldsymbol{\lambda} = (1, 1)] = 1,$$

$$\nu_{21}[(\mu, \boldsymbol{\lambda}): \mu = 5/4, \boldsymbol{\lambda} = (1,1)] = p,$$
$$\nu_{21}[(\mu, \boldsymbol{\lambda}): \mu = 5, \boldsymbol{\lambda} = (1,1)] = 1 - p.$$

Now, elementary computations show that $L_1 = (4p - 2)\log 2$, $L_2 = 0$, $\pi_1 = \pi_2 = 1/2$ so that $\sum \pi_i L_i < 0$ is equivalent to $p < 1/2$. Further, $M_{11}(s) = M_{22}(s) = 0$, $M_{12}(s) = 1$ and $M_{21}(s) = p4^s + (1-p)(1/4)^s$ so that

$$\eta(s) = [p4^s + (1-p)(1/4)^s]^{1/2} \quad \text{and} \quad \eta'(0) = (2p - 1)\log 2.$$

By Theorem 2.1, the process is certain to return to the empty state $\varnothing$ when $\sum \pi_i L_i < 0$, that is, $p < 1/2$ and has positive probability of never returning when $p > 1/2$. The condition $p < 1/2$ implies that $\eta'(0) < 0$ and so $\eta(s) = 1$ has a root

$$s_p = \frac{\log(1-p) - \log p}{\log 4} > 0$$

in this case. By Theorem 2.2, $\mathbf{E}(\tau^s)$ is finite when $0 \leq s < s_p$, and is infinite when $s > s_p$. We see that $s_p = 1$ for $p = 1/5$ and, thus, the system is transient for $p \in (1/2, 1]$, null recurrent for $p \in [1/5, 1/2)$ and positive recurrent for $p \in [0, 1/5)$, so, indeed, we have a "thick" null recurrent phase here.

In order to establish stability results for the queueing system, we need various Lyapunov functions. We can find them by studying the long term behavior of particles moving according to some random vector field models built using the parameter regeneration process and, in particular, the mean drift vector field of the process $\xi$ over service epochs.

The phenomenon of thickness of the critical region implied by the Theorem 2.2 and exhibited in the previous example is quite exceptional in the theory both of stochastic processes and of critical phenomena. This phenomenon appeared in a related model but with fixed parameters that was studied by Fayolle, Ignatyuk, Malyshev and Menshikov [4] (their results also appear in Section 5 of [3]). Their model is a random walk $\xi$ on a 2-dimensional complex, that is, a collection of quarter planes connected at their edges in some fashion with possible multiple connections and with no edge unconnected. On any quarter plane the walk drifts away from one axis and toward the other. When the boundary is reached, the quarter plane to enter next is chosen according to a transition matrix $P$. Specifically, Theorem 5.3.2 of [3] establishes recurrence of the Markov chain $\xi$ in a *fixed* environment (which satisfies our Condition E) if and only if the condition $\sum_i \pi_i L_i < 0$ holds (in fact, their result is established also for cases where the transition matrix $P$ is not irreducible). Further, Theorem 5.3.4 of [3] shows that under the same regularity conditions $\xi$ (again with fixed parameters) is positive recurrent if and only if $\eta(1) < 1$.

The new features in this paper are the following: all the parameters in the FIMM model are known constants, while here they change at each switching event according to some known distributions; moreover, the results in [4] concern only estimates of $\mathbf{E}(\tau)$, while here we have conditions for deciding whether $\mathbf{E}(\tau^s) < \infty$ or not for any $s > 0$. Finally, the phenomenon of thickness of the critical region both in [4] and in the present paper is due to some underlying randomness in the choice of the next set of parameters. In [4] this randomness is due to a Markovian choice of the next complex, analogous to our matrix $P = [P_{ij}]$: when this randomness is eliminated, the critical region shrinks to a zero measure manifold. In the present paper the thickness persists even when the Markov process induced by the matrix $P = [P_{ij}]$ becomes a deterministic routing process.

**3. The multiplicative model and the vector field model.** We first establish our results for the multiplicative model. We then define the random vector field model and show that our results for the multiplicative model enable us to describe its behavior. Using the vector field model, we can define the "random" Lyapunov function we use to establish our results for the stochastic process $(W, \xi)$.

Let $(I, X)$ denote the multiplicative model and $\sigma_A = \min\{r > 0 : X_r \leq A\}$ for any given constant $A > 0$. Part of the next result appears in [2] where conditions were found determining when the expected total time of $(I, X)$ is finite.

THEOREM 3.1. *Suppose that the multiplicative Markov chain $(I, X)$ has distributions $F_{ij}$ that satisfy Condition* E *above:*

(i) *If $\sum_1^d \pi_i L_i < 0$, then $X_r \to 0$ a.s. as $r \to \infty$, the total time $T$ of $(I, X)$ is a.s. finite and for any $s > 0$ satisfies*

$$\mathbf{E}(T^s) \begin{cases} < \infty, & \text{if } \eta(s) < 1, \\ = \infty, & \text{if } \eta(s) \geq 1. \end{cases}$$

(ii) *If $\sum_1^d \pi_i L_i = 0$, then, for any $A > 0$, $\mathbf{P}(\sigma_A < \infty \mid X_0 = x_0) = 1$ for any finite $x_0$.*

(iii) *If $\sum_1^d \pi_i L_i > 0$, then $X_r \to \infty$ a.s. as $r \to \infty$.*

REMARK. It is important to note that $\eta(0) = 1$ and $\eta$ is log-convex in $s$ so $\eta(s) < 1$ for some $s > 0$ is equivalent to $\eta'(0) = \sum \pi_i L_i < 0$. If $\sum_1^d \pi_i L_i > 0$, then $T = \infty$ a.s.

3.1. *Proofs for the multiplicative model.* We prove first the results about the limiting behavior of the process $\{X_r\}$ and then establish conditions for the existence of $\mathbf{E}(T^s)$.

PROOF OF THEOREM 3.1 (*Limiting properties of* $\{X_r\}$). The process
$$(I, \log X) \equiv \{(I_r, \log X_r) : r \geq 0\}$$
is a nonlattice random walk on $\mathcal{S}_d \times \mathbb{R}$ which is homogeneous in its second coordinate. Its jumps in this coordinate are bounded due to Condition E(ii). As we now show, its stability conditions can be demonstrated using the Lyapunov function $f(i, x) = a_i + \log x$, where the $a_i$, $i = 1, \ldots, d$, are constants to be determined. For all states $(i, x)$, let
$$\beta_i = \mathbf{E}(f(X_{r+1}) - f(X_r) \mid I_r = i, X_r = x).$$
Writing $a = (a_1, \ldots, a_d)$, $L = (L_1, \ldots, L_d)$ and $\beta = (\beta_1, \ldots, \beta_d)$, we find that $\beta = -(\mathbf{I} - P)a + L$, where $\mathbf{I}$ denotes the $d \times d$ identity matrix. The Markov chain $\{I_r\}$ is irreducible with equilibrium vector $\pi$ so it follows that the matrix $\mathbf{I} - P$ has rank $d - 1$ and that $\pi$ is orthogonal to the range of $\mathbf{I} - P$. Hence, by standard linear algebra arguments, for any $\varepsilon > 0$, there exists $a$ such that:

(i) $\beta_i \leq -\varepsilon$ for each $i$ only when $\sum_{i=1}^d \pi_i L_i < 0$;
(ii) $\beta_i = 0$ for each $i$ only when $\sum_{i=1}^d \pi_i L_i = 0$;
(iii) $\beta_i \geq \varepsilon$ for each $i$ only when $\sum_{i=1}^d \pi_i L_i > 0$.

Applying Theorem 2.1.7 of [3], we see that $\log X_r \to \pm \infty$ a.s. or, equivalently, $X_r \to \infty$ or $0$ a.s. according as to whether $\sum_{i=1}^d \pi_i L_i$ is positive or negative. If $\sum_{i=1}^d \pi_i L_i = 0$, then there exist constants $a_i$ such that $\{f(I_r, X_r)\}$ is a martingale with bounded jumps. By standard results, for example, Theorem 8.4.3 of Meyn and Tweedie [12], the process $f(I, X)$ is recurrent and as
$$\{X_r \leq A\} \supset \{f(I_r, X_r) \leq a_{\min} + \log A\},$$
the a.s. finiteness of $\sigma_A$ follows. $\square$

REMARK. A similar model, a Markov chain on the half-strip $\mathcal{S}_d \times \mathbb{Z}_+$, was studied in Section 3.1 of [3] and we borrowed the Lyapunov function $x + a_i$ from there.

We now consider the question of the *total time* of the multiplicative process. We start by restating its definition (3) in a slightly different way. For each pair $(i, j) \in \mathcal{S}_d^2$ with $i \neq j$, $F_{ij}$ is a distribution function satisfying Condition E. We construct a sequence of random variables $\alpha_1, \alpha_2, \ldots$ tied to the Markov chain $\{I_r\}$ as follows: when $I_{r-1} = i$, $I_r = j$, then $\alpha_r$ has distribution function $F_{ij}$. The $\alpha_i$ are conditionally independent given the state of $\{I_r\}$. The total time of the multiplicative model started from state $(i, 1)$ is
$$T = 1 + \alpha_1 + \alpha_1 \alpha_2 + \alpha_1 \alpha_2 \alpha_3 + \cdots = 1 + \sum_{r=1}^{\infty} \prod_{1}^{r} \alpha_n.$$

If $\sum_i \pi_i L_i \geq 0$, then the set of states $(i, x)$ with $x \geq 1$ is visited infinitely often and the total time $T$ will be infinite a.s. When $\sum_i \pi_i L_i < 0$, we can calculate which $s > 0$ have $\mathbf{E}(T^s) < \infty$.

Recall the $d \times d$ matrix $M(s) = [P_{ij}\mathbf{E}(Y_{ij}^s)]_{i,j}$, where $Y_{ij}$ has distribution function $F_{ij}$. Let $v_i(s)$ denote the right eigenvectors and $\eta_i(s)$ the corresponding eigenvalues of $M(s)$. As $\mathbf{E}(Y_{ij}^s) > 0$ for all pairs $i$, $j$ and $\{I_r\}$ is irreducible, the results of Perron–Frobenius state there is a simple eigenvalue $\eta(s) > 0$ such that $\eta(s) \geq |\eta_i(s)|$ for all other eigenvalues and its eigenvector $v(s)$ has elements $v^{(j)}(s) > 0$ and is the only eigenvector with this property.

We start by determining the asymptotic behavior of $\mathbf{E}(\alpha_1 \cdots \alpha_n)^s$.

LEMMA 3.2. *For any $s > 0$, there exist constants $C > 0$, $\gamma \in (0,1)$ such that*

$$\lim_{n \to \infty} \eta(s)^{-n}(\mathbf{E}(\alpha_1 \cdots \alpha_n)^s) = C + O(\gamma^n).$$

PROOF. Let $u$ be the row vector with $u_i = \mathbf{P}(I_0 = i)$ and $\mathbf{1}$ denote the $d$-vector of 1s. Following the standard argument for establishing the Chapman–Kolmogorov conditions, we can write

$$\mathbf{E}(\alpha_1 \cdots \alpha_n)^s = u(M(s))^n \mathbf{1}$$

and the result follows directly from this for each fixed $s$ using the result of Frobenius. □

This estimate makes it evident that the value of $\eta(s)$ is crucial to deciding whether $\mathbf{E}(T^s)$ is finite or infinite.

PROOF OF THEOREM 3.1 [*Conditions for $\mathbf{E}(T^s) = \infty$*]. The proof that $\mathbf{E}(T^s) = \infty$ when $\eta(s) > 1$ is straightforward. As the $\alpha_i$ are a.s. positive, we have for any $n$

$$\mathbf{E}(T^s) \geq \mathbf{E}(\alpha_1 \cdots \alpha_n)^s.$$

If $\eta(s) > 1$, then $\mathbf{E}(\alpha_1 \cdots \alpha_n)^s \to \infty$ as $n \to \infty$ by Lemma 3.2 and so $\mathbf{E}(T^s) = \infty$. When $\eta(s) = 1$ and $s \geq 1$, we have $T^s = (\sum_{n=1}^{\infty} \alpha_1 \cdots \alpha_n)^s \geq \sum_{n=1}^{\infty} (\alpha_1 \cdots \alpha_n)^s$ and again the result follows from Lemma 3.2.

When $\eta(s) = 1$ and $s < 1$, we can use Jensen's inequality applied to $f(x) = x^s$ which says that, for any $u$, $v$, $G$, $H > 0$ and $s < 1$ with $G \neq H$,

$$\left(\frac{uG + vH}{u+v}\right)^s > \frac{uG^s + vH^s}{u+v}$$

or

(6) $$(uG + vH)^s > (u+v)^{s-1}(uG^s + vH^s).$$

Let $T_n = 1 + \alpha_1 + \alpha_1\alpha_2 + \cdots + \alpha_1\cdots\alpha_n$ and let $t_n = \mathbf{E}(T_n^s)$. Applying (6) with $u=1$, $v=n^{-a}$, $G=T_n$ and $H=n^a\alpha_1\cdots\alpha_{n+1}$ where $a>1$ is constant, we obtain

$$t_{n+1} = \mathbf{E}(T_{n+1}^s) = \mathbf{E}((T_n + \alpha_1\cdots\alpha_{n+1})^s)$$
$$> (1+n^{-a})^{s-1}(\mathbf{E}(T_n^s) + n^{a(s-1)}\mathbf{E}(\alpha_1\cdots\alpha_{n+1})^s).$$

From Lemma 3.2, there exist constants $C>0$, $\gamma \in (0,1)$ such that the term $\mathbf{E}(\alpha_1\cdots\alpha_{n+1})^s = C + O(\gamma^n)$ and so

$$t_{n+1} > (1+n^{-a})^{s-1}(t_n + Cn^{a(s-1)}) + O(\gamma^n)$$
$$> (1+(s-1)n^{-a})(t_n + Cn^{a(s-1)}) + O(\gamma^n)$$
$$= t_n + (s-1)n^{-a}t_n + C(n^{a(s-1)} - (1-s)n^{a(s-2)}) + O(\gamma^n),$$

where $(1+n^{-a})^{s-1} > 1+(s-1)n^{-a}$ as $s<1$. Choose $a = (1-s/2)/(1-s) > 1$ and for some constant $\beta \in (0, C/s)$, pick $n_0 \equiv n_0(C,s)$ such that, for $n \geq n_0$,

$$(1-s)(\beta n^{-a+s/2} + Cn^{a(s-2)}) + O(\gamma^n) < (C - \beta s/2)n^{-1+s/2}.$$

This is possible as $-a + s/2 < -1 + s/2$ and $a(s-2) < a(s-1) = -1 + s/2$. As $t_n > t_1 = \mathbf{E}(\alpha_1^s) > 0$ for all $n$, we can choose $\beta \in (0, t_{n_0}n_0^{-s/2})$. Hence, $t_n > \beta n^{s/2}$ for some $n \geq n_0$ and for such $n$,

$$t_{n+1} > \beta n^{s/2} + Cn^{-1+s/2} - (1-s)(\beta n^{-a+s/2} + Cn^{a(s-2)}) + O(\gamma^n)$$
$$> \beta\left(n^{s/2} + \frac{s}{2}n^{-1+s/2}\right) > \beta(n+1)^{s/2},$$

so, by induction, $t_n > \beta n^{s/2}$ for all $n \geq n_0$. Finally, as $T^s > T_n^s$ for all $n$, $\mathbf{E}(T^s) = \infty$.

*Conditions for $\mathbf{E}(T^s) < \infty$.* Now we assume $\eta(s) < 1$. If $s \leq 1$, then for, any $n$,

$$T^s = \left(\sum_{n=1}^{\infty} \alpha_1\cdots\alpha_n\right)^s \leq \sum_{n=1}^{\infty}(\alpha_1\cdots\alpha_n)^s$$

and $E(T^s) < \infty$ by comparison with the tail sum $\sum_{n \geq n_0} \eta(s)^n$.

It only remains to sort out the case when $s > 1$. Again, applying Jensen's inequality to the function $f(x) = x^s$, we have, for any real numbers $u$, $v$, $G$, $H > 0$ and $s > 1$ with $G \neq H$,

$$\left(\frac{uG + vH}{u+v}\right)^s < \frac{uG^s + vH^s}{u+v}$$

or

(7) $$(uG + vH)^s < (u+v)^{s-1}(uG^s + vH^s).$$

Fix any $\beta \in ((\eta(s))^{1/s}, 1)$, and apply (7) with $u = 1$, $v = \beta^{n+1}$, $G = T_n$, $H = \beta^{-(n+1)}\alpha_1 \cdots \alpha_{n+1}$ to obtain the inequality

$$
\begin{aligned}
(8) \quad \mathbf{E}(T_{n+1}^s) &= \mathbf{E}(T_n + \beta^{n+1} \times \beta^{-(n+1)}\alpha_1 \cdots \alpha_{n+1})^s \\
&\leq (1 + \beta^{n+1})^{s-1}(\mathbf{E}(T_n^s) + \beta^{n+1}\beta^{-s(n+1)}\mathbf{E}(\alpha_1 \cdots \alpha_{n+1})^s).
\end{aligned}
$$

From the definition of $\eta(s)$ and, choice of $\beta$, it follows that

$$\lim_{n \to \infty} \beta^{-s(n+1)}\mathbf{E}(\alpha_1 \cdots \alpha_{n+1})^s = 0.$$

Again, as $t_n \geq t_1 = \mathbf{E}(\alpha_1^s) > 0$ for any $n$, it follows that there exists $n_0$ such that, for all $n \geq n_0$,

$$\beta^{-s(n+1)}\mathbf{E}(\alpha_1 \cdots \alpha_{n+1})^s \leq t_n.$$

Now (8) implies that, for all $n \geq n_0$,

$$t_{n+1} \leq (1 + \beta^{n+1})^{s-1}(t_n + \beta^{n+1}t_n) = (1 + \beta^{n+1})^s t_n.$$

Thus, for any $n$,

$$\mathbf{E}(T_n^s) \leq t_{n_0} \prod_{m=n_0}^{\infty} (1 + \beta^{m+1})^s < \infty.$$

That $\mathbf{E}(T^s) < \infty$ when $\eta(s) < 1$ and $s > 1$ follows from this bound and the monotone convergence theorem. $\square$

3.2. *Definition and properties of the random vector field model.* We now define a random vector field model $\{(S(\omega, t), V(\omega, t)) : \omega \in \Omega_S, t \geq 0\}$ using the parameter process $W$. It is motivated by the one-step mean drift vector field defined in (1), but we define a more general version so we can consider perturbations of this natural field. This type of model is simple and standard, but we describe it carefully here as we use it to define the *random* Lyapunov functions we need to establish the properties of the random walk process $\xi$.

There are no switching times in the vector field model. Consider a particle that moves around on $\mathcal{S}_d^2 \times \mathbb{R}_+^2$ with velocity field $\mathbf{D} = \{\mathbf{D}(\omega, r) \in \mathbb{R}_- \times \mathbb{R}_+ : \omega \in \Omega_S, r \in \mathbb{N}\}$, where the $\mathbf{D}(\omega, r) = (-d_1(\omega_r), d_2(\omega_r))$ satisfy conditions equivalent to Condition E which we now specify.

CONDITION E′. As before, let $F_{ij}$ denote the distribution function of the slope $Y(\omega, r) = d_2(\omega_r)/d_1(\omega_r)$ and assume that there exists some constants $m_0$, $M_0$ such that

$$\nu_{ij}[d_1 > m_0, \|\mathbf{D}\| < M_0] = 1 \quad \text{and that} \quad F_{ij}(1/M_0) = 0, \ F_{ij}(M_0) = 1$$

for each pair $i$, $j$.

For our regenerated parameter process $(W, \xi)$, the *natural* vector field satisfies $d_1 = 1 - \lambda_1/\mu$, $d_2 = \lambda_2/\mu$ and we denote it by $\mathbf{D}^0$. In this section we are considering $d_1$ and $d_2$ to be more general functions of the parameters.

The process $(S, V)$ on the state space $\Omega_e \times \mathbb{N} \times \mathbb{R}_+^2$ has components

$$S(\omega, t) = (i_r, i_{r+1}, G_r, \boldsymbol{\lambda}_r, r),$$

the parameter set at time $t$ including the current epoch number and

$$V(\omega, t) = (v_1(\omega, t), v_2(\omega, t)),$$

the position of the particle on the quarter plane indicated by the parameter $i_r$ at time $t$. For each parameter sequence $\omega$ and initial point $(x_0, 0)$ with $x_0 > 0$ and $t_0 = 0$, set $S(\omega, t) = S(\omega, 0) = \omega_1 = (i_1, i_2, G_1, \boldsymbol{\lambda}_1, 1)$ and

$$V(\omega, t) = (x_0, 0) + t\mathbf{D}(\omega, 1)$$

until the particle reaches the boundary of quarter plane $(i_1, i_2)$ which defines the end of the current epoch. This occurs at time

$$t_1 \equiv t_1(\omega, x_0) = \frac{x_0}{d_1(\omega_1)}.$$

Switching times are assumed to be 0 in this model so at $t = t_1$ set

$$S(\omega, t_1) = \omega_2 \quad \text{and} \quad V(\omega, t_1) = \left(d_2(\omega_1)\frac{x_0}{d_1(\omega_1)}, 0\right) = (x_0 Y(\omega, 1), 0),$$

where $Y(\omega, r) = d_2(\omega_r)/d_1(\omega_r)$. The components of $V$ are transposed at $t_1$ as the particle moves from plane $(i_1, i_2)$ to plane $(i_2, i_3)$ at that moment, corresponding to the server switching stations. For epochs $r = 2, 3, \ldots$, define their end times by

$$(9) \qquad t_r \equiv t_r(\omega, x_0) = t_{r-1} + \frac{v_1(\omega, t_{r-1})}{d_1(\omega_r)}.$$

Now, for $t \in [t_{r-1}, t_r)$, we set $S(\omega, t) = \omega_r$ and

$$V(\omega, t) = V(\omega, t_{r-1}) + (t - t_{r-1})\mathbf{D}(\omega, r)$$

so that $v_2(\omega, t) = d_2(\omega_r)(t - t_{r-1})$ and then

$$V(\omega, t_r) = (d_2(\omega_r)(t_r - t_{r-1}), 0).$$

It is clear from these definitions that, for every parameter sequence $\omega$, the sequence $\{t_r\}$ is increasing and so $t_\varnothing(\omega, x_0) = \lim_{r \to \infty} t_r(\omega, x_0)$ exists. Under conditions such that $t_\varnothing$ is bounded a.s., it is the time taken for the particle to reach the empty state $\varnothing$. We now determine conditions which ensure $\|V(\omega, t)\| \to 0$ and conditions under which $\|V(\omega, t)\| \to \infty$ as $t \to \infty$ and then consider the process time.

The sequence $v_1(\omega, t_r)$, $r = 1, 2, \ldots$, describes the distance of the particle from the origin at the starts of epochs $r + 1$. For each $r = 1, 2, \ldots$,

$$v_1(\omega, t_r) = d_2(\omega_r)(t_r - t_{r-1}) = d_2(\omega_r)\frac{v_1(\omega, t_{r-1})}{d_1(\omega_r)} = v_1(\omega, t_{r-1})Y(\omega, r)$$

and, hence,

$$(10) \qquad v_1(\omega, t_r) = x_0 \prod_{n=1}^{r} Y(\omega, n).$$

We re-iterate at this point that this vector field model is a random process driven by the parameter regeneration process including the Markov chain $\{I_r\}$ for the server location. Recall that when $I_n = i$, $I_{n+1} = j \neq i$, then $Y(\omega, n) = d_2(\omega_n)/d_1(\omega_n)) \sim Y_{ij}$ with distribution $F_{ij}$; $\pi = (\pi_1, \pi_2, \ldots, \pi_d)$ denotes the equilibrium distribution of $\{I_r\}$; $L_i = \sum_{j=1}^{d} P_{ij}\mathbf{E}(\log Y_{ij})$ for each $i$. For $C > 0$, let $t_C = \min\{t \geq 0 : v_1(\omega, t) + v_2(\omega, t) \leq C\}$.

LEMMA 3.3. *For the random vector field model* $(S, V)$:

(i) *if* $\sum_{i=1}^{d} \pi_i L_i < 0$, *then* $v_1(\omega, t_r) \to 0$ *a.s. and*

$$t_\varnothing(\omega, x_0) = \lim_{r \to \infty} t_r(\omega, x_0) < \infty \qquad a.s.,$$

(ii) *if* $\sum_{i=1}^{d} \pi_i L_i = 0$, *then* $\mathbf{P}(t_C < \infty) = 1$,
(iii) *if* $\sum_{i=1}^{d} \pi_i L_i > 0$, *then* $v_1(\omega, t_r) \to \infty$ *a.s. as* $r \to \infty$.

PROOF. The process $\{(I_r, v_1(\omega, t_r))\}_r$ is a multiplicative Markov process, so parts (i) and (iii) follow immediately from Theorem 3.1(i) and (iii). Theorem 3.1(ii) implies that $\min\{r : v_1(\omega, t_r) < C\}$ is finite $\omega$-a.s. From equation (9) epoch $r$ lasts for time $t_r - t_{r-1} = v_1(\omega, t_{r-1})/d_1(\omega_r)$ which is finite as $d_1(\omega_r) \geq m_0$ and, hence, $t_C < \infty$ a.s. $\square$

REMARK. As $t_r(\omega, x_0) = x_0 t_r(\omega, 1)$ for all $r$ and all $x_0 > 0$, it follows that if $t_\varnothing(\omega, 1) < \infty$, then $t_\varnothing(\omega, x_0) < \infty$ for all $x_0 > 0$.

Next we recall the Perron–Frobenius eigenvalue $\eta(s)$ of the matrix $M(s)$ defined in (4).

LEMMA 3.4. *Suppose* $\sum_{i=1}^{d} \pi_i L_i < 0$. *For* $s > 0$, *if* $\eta(s) < 1$, *then* $\mathbf{E}(t_\varnothing^s) < \infty$ *while if* $\eta(s) \geq 1$, *then* $\mathbf{E}(t_\varnothing^s) = \infty$.

PROOF. Combining equations (9) and (10), we obtain

$$t_r - t_{r-1} = \frac{v_1(\omega, t_{r-1})}{d_1(\omega_r)} = \frac{v_1(\omega, t_r)}{d_2(\omega_r)} = \frac{x_0}{d_2(\omega_r)} \prod_{n=1}^{r} Y(\omega, n)$$

for $r = 1, 2, \ldots$ and so

$$t_\varnothing = \sum_{r=1}^\infty t_r - t_{r-1} = x_0 \sum_{r=1}^\infty \frac{1}{d_2(\omega_r)} \prod_{n=1}^r Y(\omega, n).$$

Under Condition $\mathrm{E}'$ for the vector field model, there exists a constant $K > 0$ (uniformly in $\omega$) such that $1/K < 1/d_2(\omega_r) < K$ a.s. for all epochs $r$. It follows that

$$\frac{x_0}{K} \sum_{r=1}^\infty \prod_1^r Y(\omega, n) < t_\varnothing(\omega, x_0) < K x_0 \sum_{r=1}^\infty \prod_1^r Y(\omega, n).$$

The expression $\sum_r \prod_n Y(\omega, n)$ is the total time for the multiplicative random walk associated with the parameter regeneration process that drives this vector field model. From Theorem 3.1(i), it follows that $\mathbf{E}(t_\varnothing^s) < \infty$ is equivalent to $\eta(s) < 1$. $\square$

The assumptions in Condition $\mathrm{E}'$ ensure this construction can be carried out for indefinitely many epochs for all parameter sequences from any initial particle position $\mathbf{x}$ with $x > 0$ in any epoch $r$. To indicate this, we denote the switching instants defined above by $t_n(\omega, r, \mathbf{x})$ for $n \geq r$. Under the condition $\sum_{i=1}^d \pi_i L_i < 0$, it follows that $t_\varnothing(\omega, r, \mathbf{x}) = \lim_{n \to \infty} t_n(\omega, r, \mathbf{x}) < \infty$ and is the time taken by the particle to reach 0 from initial point $\mathbf{x}$ in epoch $r$. In this case, for any constant $A > 0$,

(11) $\quad t_A(\omega, r, \mathbf{x}) = \min\{t : V(t) \in \mathbf{B}_A, S(0) = \omega_r, V(0) = \mathbf{x}\} < \infty \qquad \text{a.s.},$

where $\mathbf{B}_A = \{\mathbf{x} \in \mathbb{R}_+^2 : x + y \leq A\}$.

**4. Lyapunov functions and their properties.** Now we return to studying $(W, \xi)$, the random walk together with its regenerated parameter process. We use the vector field model to construct Lyapunov functions along each parameter sequence $\omega$ and then take expectations over the parameter process to establish Theorems 2.1 and 2.2. We outline the simplest case of our scheme now to help explain this general argument. Let $\mathbf{E}_\omega$ denote expectation for the random walk conditioned on the parameter sequence $\omega$, let $(r_n, \xi_n)$ denote the epoch and queue lengths at time $n$ and let $\tau_A = \min\{n : \xi_n \in \mathbf{B}_A\}$. For the natural vector field $\mathbf{D}^0(\omega, r) \equiv (-d_1^0(\omega, r), d_2^0(\omega, r)) = (\lambda_{r,1}/\mu_r - 1, \lambda_{r,2}/\mu_r)$ consider the function $t_A$ defined in (11). For $\omega$ such that $t_A$ is finite and for $\mathbf{x} \notin \mathbf{B}_A$,

(12) $\qquad \mathbf{E}_\omega(t_A(\omega, r, \xi_{n+1}) - t_A(\omega, r, \xi_n) \mid \xi_n = \mathbf{x}) = -1.$

In the case where switching happens instantaneously, Foster's criterion (see Theorems 2.1.1 and 2.2.3 of [3]) implies that

$$\mathbf{E}_\omega(\tau_A \mid \xi_0 = \mathbf{x}_0) < t_A(\omega, 1, \mathbf{x}_0) < \infty$$

and it follows, by taking expectation over $\Omega_S$ and using Lemma 3.4 with $s = 1$, that $\mathbf{E}(\tau_A \,|\, \xi_0 = \mathbf{x}_0) < \infty$ whenever $\eta(1) < 1$.

We see from this that, during epoch $r$, $t_A(\omega, r, \cdot)$ based on $\mathbf{D}^0$ is a random ($\omega$-dependent) Lyapunov function for establishing stability ($\mathbf{E}(\tau_A) < \infty$) of the queue length process $\xi$. Our plan is to show that when $\mathbf{E}(t_A^s) < \infty$ (conditions for this are given in Lemma 3.4), then $\mathbf{E}(\tau_A^s) < \infty$. However, a stronger result than Foster's is necessary to decide whether $\mathbf{E}(\tau_A^s) < \infty$ for given $s > 1$ (see Theorem 5.3 below). This result requires stronger bounds, than (12). To get these stronger bounds we actually have to use $t_A$ from a vector field which is a modification of $\mathbf{D}^0$. It is also necessary to handle nonzero switching times.

To establish which moments of $\tau_A$ are infinite, a different approach is needed. We show that for almost all parameter sequences $\omega$ there is some probability, $\beta > 0$, say, that the random walk takes time comparable with the vector field to reach $\mathbf{B}_A$. This $\beta$ is uniform over $\omega$, so we can take expectation over $\Omega_S$ to show that $\mathbf{E}(t_A^s) = \infty$ implies $\mathbf{E}(\tau_A^s) = \infty$.

4.1. *Modified vector fields and Lyapunov functions.* Consider the set of particle trajectories, under a vector field $\mathbf{D}$, from all possible initial points on the $\mathbf{e}_1$ axis in epoch 1. For any parameter sequence $\omega$, any point $(r, \mathbf{x}) \in \mathbb{N} \times \mathbb{R}_+^2$ lies on a unique such trajectory.

DEFINITION 4.1. For a vector field $\mathbf{D}$ satisfying Condition E′, we define two distinct Lyapunov functions. For each parameter sequence $\omega$, define

$$g : \Omega_S \times \mathbb{N} \times \mathbb{R}_+^2 \to \mathbb{R}_+ \qquad \text{by } g(\omega, r, \mathbf{x}) = \gamma,$$

where $\gamma$ is such that the trajectory under $\mathbf{D}(\omega, \cdot)$, started from $(1, \gamma \mathbf{e}_1)$, passes through the point $(r, \mathbf{x})$, that is, $g$ is constant along particle trajectories under the vector field $\mathbf{D}$. On

$$F = \{\omega \in \Omega_S : t_\varnothing(\omega, 1, \mathbf{e}_1) < \infty\},$$

the set of parameter sequences for which $t_\varnothing$ is everywhere finite, define $f : F \times \mathbb{N} \times \mathbb{R}_+^2 \to \mathbb{R}_+$ by

$$f(\omega, r, \mathbf{x}) = t_A(\omega, r, \mathbf{x}) \qquad \text{where } A > 0 \text{ is a given constant.}$$

We will need to consider expectations $\mathbf{E}_\omega$ of $f$ and $g$ as functions of the random walk $\xi$ for given $\omega$ and, subsequently, with respect to the measure on the parameter sequences.

REMARK. For any given $\omega$, the functions $f$ and $g$ are linear within epochs and we consider their properties for general vector fields $\mathbf{D}$. We will call $f$ the *remaining life*, as it is the time for the particle to reach $\mathbf{B}_A = \{\mathbf{x} \in \mathbb{R}_+^2 : x + y \leq A\}$ from initial point $\mathbf{x}$ in epoch $r$ along the trajectory defined by the vector field $\mathbf{D}$ for the parameter sequence $\omega \in F$. Similarly, we will call $g$ the *initial value* function and observe that it is everywhere finite for a.s. all $\omega$ for any $\mathbf{D}$ satisfying Condition E$'$.

In particular, we will use $f$ for establishing that specified moments of the random time $\tau_A = \min\{n : \xi_n \in \mathbf{B}_A\}$ exist, $g$ for establishing instability and $f$ and $g$ together for showing nonexistence of certain moments of $\tau_A$. We now describe some of their properties.

LEMMA 4.2. *If vector field $\mathbf{D}$ has $\sum_i \pi_i L_i < 0$, then $P_\Omega(F) = 1$. For fixed $\omega$, the functions $f$ (for $\omega \in F$) and $g$ are linear in $\mathbf{x}$ for fixed $r$ and continuous at changes of epoch $r$. Writing:*

(i) $f(\omega, r, \mathbf{x}) = a_r x_1 + b_r x_2 + f_0$ *for $\omega \in F$ ($a_r$, $b_r$ and $f_0$ are functions of $\omega$), we have:*

- $a_r = b_{r-1}$ *for each $r$ and there exists $M > 0$, uniform in $\omega$ and $r$, such that $b_r/M \leq a_r \leq M b_r$,*
- $\phi_r \equiv \max\{a_r, b_r\} \geq 1/M_0$ *(the constant in Condition E$'$),*
- $f(\omega, r, \mathbf{x}) \geq (x_1 + x_2)/M_0$;

(ii) $g(\omega, r, \mathbf{x}) = c_r x_1 + c_{r+1} x_2$, *we have:*

- $c_{r+1} = c_r / Y(\omega, r)$ *for each $r$ and, hence, $c_r / c_{r+1} = Y(\omega, r) \in (1/M_0, M_0)$ for all epochs $r$,*
- $c_r \prod_{n=1}^{r-1} Y(\omega, n) = c_1$ *for all $r$.*

PROOF. $P_\Omega(F) = 1$ is the content of Lemma 3.3(i). That $f$, for $\omega \in F$, and $g$ are linear in $\mathbf{x}$ for fixed $r$ and continuous at changes of epoch $r$ follows from the vector field $\mathbf{D}$ being constant during epochs.

The properties of $f$ for $\omega \in F$ follow from Condition E$'$ which constrains $d_1$ and $d_2$ and the relations $a_r(x_1 - d_1) + b_r(x_2 + d_2) = a_r x_1 + b_r x_2 - 1$ or, equivalently, $1 + b_r d_2(r) = a_r d_1(r)$ for each $r$. The lower bound on $f$ is immediate. The relations satisfied by the $c_r$ follow similarly. □

Unfortunately, $a_r$ and $b_r$, the $\mathbf{x}$-coefficients of $f$ within an epoch, can be arbitrarily large, so the level curves of $f(\omega, r, \cdot)$ (constructed under $\mathbf{D}^0$) in the $\mathbf{x}$-plane can make an arbitrarily small angle with the natural drift direction $\mathbf{D}^0(\omega, r)$ of $\{\xi_n\}$. This prevents us making the estimates needed to establish which hitting time moments exist if we construct $f$ with respect to field $\mathbf{D}^0$.

*Modified vector fields.* To resolve this problem, we construct modified vector fields (small perturbations of the natural vector field $\mathbf{D}^0$) for which $f$ and $g$ permit some useful estimations. We do this in different ways for showing whether $\mathbf{E}(\tau_A^s) < \infty$ or not.

To show that $\mathbf{E}(\tau_A^s) < \infty$, we want a vector field $\mathbf{D}$ such that its remaining life $f$ is greater than that of $\mathbf{D}^0$ but is still finite a.s. and for which moment $s$ of $f$ exists. When we want to show that $\mathbf{E}(\tau_A^s) = \infty$ for some $s$, we use a different modification of $\mathbf{D}^0$ and the initial value function $g$.

DEFINITION 4.3. As in equation (2), we write
$$\theta_0(\omega, r) = \arctan[\lambda_2(\omega, r)/(\mu(\omega, r) - \lambda_1(\omega, r))]$$
for the minor angle between $\mathbf{D}^0(\omega, r)$ and $-\mathbf{e}_1$. We define

(13)
$$\mathbf{D}^+(\omega, r) \equiv (-d_1^+(\omega_r), d_2^+(\omega_r))$$
$$= \|\mathbf{D}^0(\omega, r)\|(-\cos(\theta_0(\omega, r) + \theta'), \sin(\theta_0(\omega, r) + \theta'))$$

for some $\theta' > 0$ such that $\theta_0(\omega, r) + \theta' < \pi/2$ for all pairs $\omega$, $r$. We will say that such a $\mathbf{D}^+$ is *uniformly above* $\mathbf{D}^0$. We further define

(14)
$$\mathbf{D}^-(\omega, r) \equiv (-d_1(\omega_r), d_2(\omega_r))$$
$$= \|\mathbf{D}^0(\omega, r)\|(-\cos(\theta_0(\omega, r) - \theta'), \sin(\theta_0(\omega, r) - \theta'))$$

for some sufficiently small $\theta' > 0$. In this case we say that $\mathbf{D}^-$ is *uniformly below* $\mathbf{D}^0$.

This definition is possible, as, by Condition E,
$$\theta_0(\omega, r) \in (\arctan m_0/M_0, \arctan M_0/m_0) \subset (0, \pi/2).$$

Clearly, there are $\omega$ for which the remaining life under $\mathbf{D}^0$ is finite, but under $\mathbf{D}^+$ it becomes infinite. In fact, the set of these has zero probability for some suitably small $\theta'$ as we now show. Let $M^+(s)$, $M^0(s)$ and $M^-(s)$ denote the matrix $[P_{ij}\mathbf{E}(Y_{ij}^s)]$ under the vector fields $\mathbf{D}^+$, $\mathbf{D}^0$ and $\mathbf{D}^-$ respectively and let $\eta^+(s)$, $\eta^0(s)$ and $\eta^-(s)$ denote the Perron–Frobenius eigenvalues of these matrices. Also, for each station $i$, let $L_i^0 = \sum_{j=1}^d P_{ij}\mathbf{E}(\log Y_{ij})$ under $\mathbf{D}^0$ and $L_i^-$ the same quantities under $\mathbf{D}^-$.

LEMMA 4.4. *For any given $s > 0$, we can choose $\theta' > 0$ small enough that:*

  (i) *if $\eta^0(s) < 1$, then $\eta^+(s) < 1$;*
 (ii) *if $\eta^0(s) > 1$, then $\eta^-(s) > 1$;*
(iii) *if $\sum_i \pi_i L_i^0 > 0$, then $\sum_i \pi_i L_i^- > 0$*

*($\theta'$ does not depend upon $\omega$).*

REMARK. When $\sum_i \pi_i L_i^0 < 0$ and $\theta'$ is chosen to satisfy condition (i) or (ii), we say that $\mathbf{D}^+$ or $\mathbf{D}^-$ respectively is *s-neutral*.

PROOF OF LEMMA 4.4. For the vector field $\mathbf{D}^+$, we have $Y_{ij} = d_2^+/d_1^+ = \tan(\theta_0 + \theta')$ when $\omega_r = (i, j, G, \boldsymbol{\lambda})$. As the function tan is continuous on $(0, \pi/2)$ and $\eta^+(s)$ is a continuous function of $M^+(s)$, it follows that $\eta^+(s)$ is continuous in $\theta'$. Part (i) follows as Condition E bounds the rate of change of $\eta^+(s)$ in $\theta'$ and parts (ii) and (iii) follow from a similar argument. $\square$

The advantage of working with $\mathbf{D}^+$ rather than $\mathbf{D}^0$ is that the expected one-step jumps of the process $\{f(\omega, r_n, \xi_n)\}_n$ are of size comparable to $\phi_r$ where, as before, $f(\omega, r, \mathbf{x}) = a_r x + b_r y + f_0$ and $\phi_r = \max\{a_r, b_r\}$. Let $\Delta f_n = f(\omega, r_{n+1}, \xi_{n+1}) - f(\omega, r_n, \xi_n)$ and $\Delta g_n = g(\omega, r_{n+1}, \xi_{n+1}) - g(\omega, r_n, \xi_n)$.

LEMMA 4.5. *Suppose that $\mathbf{D}^+$ lies uniformly above $\mathbf{D}^0$ and is s-neutral for some $s > 0$. Let $f$ denote the remaining life of $\mathbf{D}^+$. Then $F^+ = \{\omega \in \Omega_S : f(\omega, 1, \mathbf{e}_1) < \infty\}$ satisfies $P_\Omega(F^+) = 1$ and for each $\omega \in F^+$ and any service epoch $r$,*

$$\mathbf{E}(f(\omega, r, \mathbf{x})^s) < \infty$$

*and there exists $\varepsilon_1 > 0$ such that*

$$\mathbf{E}_\omega(\Delta f_n \mid r_n = r, \xi_n = \mathbf{x}) = -a_r d_1^0(\omega, r) + b_r d_2^0(\omega, r) \leq -\varepsilon_1 \phi_r,$$

*where $\phi_r \geq 1/M_0$. Similarly, the initial value function $g$ of a vector field $\mathbf{D}^-$ which is uniformly below $\mathbf{D}^0$ satisfies*

$$\mathbf{E}_\omega(\Delta g_n \mid r_n = r, \xi_n = \mathbf{x}) = -c_r d_1^0(\omega, r) + c_{r+1} d_2^0(\omega, r) \geq \varepsilon_2 c_r$$

*for some $\varepsilon_2 > 0$ which does not depend upon $r$ and there exists a constant $\kappa > 0$ such that*

$$\mathbf{E}_\omega(\Delta g_n^2 \mid r_n = r, \xi_n = \mathbf{x}) \leq c_r^2 \kappa.$$

PROOF. $P(F^+) = 1$ and $\mathbf{E}(f(\omega, r, \mathbf{x})^s) < \infty$ follow immediately from Lemma 3.4 applied to the vector field $\mathbf{D}^+$.

Let $\beta$, $\beta^+$ denote the angles between the vector $(a_r, b_r)$ and $\mathbf{D}^0$, $\mathbf{D}^+$ respectively. By comparison with (12), we see that $\beta > \pi/2$ and, by construction, there exists $\theta' > 0$ such that $\beta^+ = \beta + \theta'$. The inequality for $\mathbf{E}(\Delta f_n)$ follows from the comparability of $a_r$ and $b_r$, $\phi_r \geq 1/M_0$ (from Lemma 4.2), $\|\mathbf{D}^0\| \geq m_0$ and the cosine rule.

The inequality for $\mathbf{E}(\Delta g_n)$ follows from a very similar argument, while the existence of the constant $\kappa$ for bounding $\mathbf{E}(\Delta g_n^2)$ follows from the restriction on the rates $\boldsymbol{\lambda}$ imposed by Condition E and the assumption that service time distributions $G$ are selected from a set $\Gamma$ with uniformly bounded second moments. $\square$

**5. Proofs for the polling system.** We are now in a position to establish Theorems 2.1 and 2.2. Our treatment of the vector field model in the previous sections ignored switching times, but we must now consider them.

5.1. *Technique for handling switching times.* Our main method is to observe the process not at all transitions, but at carefully calculated sampling times. Theorem 5.1 below states the semi-martingale result from Fayolle, Menshikov and Malyshev [3] which we will use to establish the recurrence part of Theorem 2.1.

Let $\{Z_t, t \geq 0\}$ be a sequence of nonnegative random variables with $Z_0$ constant and $Z_t$ measurable with respect to a filtration $\mathcal{F} = \{\mathcal{F}_t, t \geq 0\}$. Let $\{N_n, n \geq 0\}$ be a strictly increasing sequence of stopping times adapted to $\mathcal{F}$ such that $N_0 = 0$ and $|N_{n+1} - N_n|$ is bounded for all $n \geq 0$. Define a sampled process $\{X_n, n \geq 0\}$ by $X_0 = Z_0$ and $X_n = Z_{N_n}$ for $n \geq 1$. Also, for constant $C > 0$, let $\zeta_C = \min\{t \geq 1 : Z_t \leq C\}$ and $\sigma = \min\{n \geq 1 : X_n \leq C\}$. Finally, let $\{X_{t \wedge \sigma}\}$ denote the sequence $\{X_t\}$ stopped at $\sigma$. With this notation, we now state Theorem 2.1.2 of [3].

THEOREM 5.1. *If $Z_0 > C$ and for some $\varepsilon > 0$ and all $n \geq 0$*

$$\mathbf{E}(X_{(n+1) \wedge \sigma} \mid \mathcal{F}_{N_{n \wedge \sigma}}) \leq X_{n \wedge \sigma} - \varepsilon \mathbf{E}(N_{(n+1) \wedge \sigma} - N_{n \wedge \sigma} \mid \mathcal{F}_{N_{n \wedge \sigma}}) \qquad a.s.,$$

*then $\mathbf{E}(\zeta_C) \leq Z_0/\varepsilon < \infty$.*

To use this result, we will calculate in the proof of Theorem 2.1 a constant $N$ such that the inequality of this theorem is satisfied $N$ steps after any switching time.

5.2. *A key supermartingale for the instability results.* We now show that the level curves of the initial value function $g$ of $\mathbf{D}^-$ form a useful barrier for the random walk $\xi$. For each $\omega \in \Omega$, $K > 0$, let

(15) $$\mathcal{G}_K(\omega) = \{(r, \mathbf{x}) \in \mathbb{N} \times \mathbb{R}_+^2 : g(\omega, r, \mathbf{x}) \leq K\}$$

denote the set of states below the trajectory of $\mathbf{D}^-$ with initial value $K$. Let

$$\rho_K(\omega) = \min\{n : (r_n, \xi_n) \in \mathcal{G}_K(\omega)\}$$

denote the hitting time of the random walk to $\mathcal{G}_K$ ($r_n$ denotes the epoch in which step $n$ of $\xi$ occurs). We recall that $\mathbf{B}_A = \{\mathbf{x} \in \mathbb{R}_+^2 : x + y \leq A\}$.

LEMMA 5.2. *Let $g$ denote the initial value function of vector field $\mathbf{D}^-$ and let $\tau_A$ denote the hitting time of $\mathbf{B}_A$ by $\xi$. Let $Z_n = g(\omega, r_n, \xi_n)$ for $n \in \mathbb{N}$ and suppose that $Z_0 = z_0$ for some constant $z_0 > A$. Define $U_n = 1/Z_{n \wedge \tau_A}$* for $n = 1, 2, \ldots$. If $A$ is sufficiently large, then, for every epoch $r$ and every parameter sequence $\omega$,

$$\mathbf{E}_\omega(U_{n+1} - U_n \mid r_n = r, \xi_n = \mathbf{x}) < 0$$

and, hence, $U_\infty$ exists a.s.

PROOF. Taylor's theorem applied to the function $(1+x)^{-1}$ provides the inequality

$$(1+x)^{-1} - 1 < -x + x^2(1-t_0)^{-3} \qquad \text{at any } x > -t_0 > -1.$$

Write $Z_n = z$, $r_n = r$ and consider $\xi_n = \mathbf{x} \notin \mathbf{B}_A$. From Lemma 4.2(ii), we have $z \geq c_r A/M_0$ and either $\Delta g_n = Z_{n+1} - Z_n \geq 0$ or $0 < -\Delta g_n \leq c_r$. In either case $\Delta g_n/z \geq -M_0/A > -1$ for $A > M_0$. It now follows from Lemma 4.5 and the assumption $c_r/z \leq M_0/A$ that

$$\begin{aligned}
\mathbf{E}_\omega(U_{n+1} - U_n \mid \xi_n = \mathbf{x}) &= \frac{1}{z}\left(\mathbf{E}_\omega\left(\frac{1}{1+\Delta g_n/z} \,\Big|\, \xi_n = \mathbf{x}\right) - 1\right) \\
&< \frac{1}{z}\mathbf{E}_\omega\left(-\frac{\Delta g_n}{z} + (1-M_0/A)^{-3}\frac{\Delta g_n^2}{z^2} \,\Big|\, \xi_n = \mathbf{x}\right) \\
&< -\frac{\varepsilon_2 M_0}{zA} + (1-M_0/A)^{-3}\frac{\kappa M_0^2}{zA^2}.
\end{aligned}$$

The right-hand side can be made negative by choosing $A$ sufficiently large. Therefore, along any given parameter sequence, the process $(U_n)$ is a bounded super-martingale which guarantees the existence of $U_\infty$. $\square$

5.3. *Proof of Theorem* 2.1. (i) We wish to show that $\tau$, the return time to the empty state $\varnothing$, is a.s. finite under the condition $\sum_i \pi_i L_i < 0$. Consider a vector field $\mathbf{D}^+$ which is uniformly above $\mathbf{D}^0$ and $s$-neutral for some $s > 0$. As before, we denote the remaining life of $\mathbf{D}^+$ by $f$ and let $F$ denote the set of parameter sequences where $f$ is a.s. finite for all points $(r, \mathbf{x})$. For each $\omega \in F$, Lemma 4.5 implies

$$\mathbf{E}_\omega(f(\omega, r, \xi_{n+1}) - f(\omega, r, \xi_n) \mid \xi_n = \mathbf{x}) \leq -\varepsilon_1 \phi_r \leq -\varepsilon_1/M_0$$

during service epochs. During switching times, the queue lengths have drift upward so this inequality does not hold. Let $T_{r-1}$ denote the time at which the event which concludes epoch $r - 1$ occurs. Then

$$\begin{aligned}
&\mathbf{E}_\omega(f(\omega, r, \xi_{n+1}) - f(\omega, r, \xi_n) \mid \xi_n = \mathbf{x}, T_{r-1} = n, \omega_r = (i, j, G, \boldsymbol{\lambda})) \\
&= (a_r\lambda_1 + b_r\lambda_2)m_{ij} > 0,
\end{aligned}$$

where $m_{ij}$ is the mean switching time from station $i$ to $j$. To employ Theorem 5.1, we choose $\varepsilon_1 > 0$ according to Lemma 4.5, $\varepsilon' \in (0, \varepsilon_1)$ and an integer $N \geq 2M_0 m^*/(\varepsilon_1 - \varepsilon')$, where $m^* = \max_{ij} m_{ij}$. As long as $\mathbf{x} = (x,y)$ has $x > N$,

$$\begin{aligned}
&\mathbf{E}_\omega(f(\omega, r, \xi_{n+N}) - f(\omega, r, \xi_n) \mid \xi_n = \mathbf{x}, T_{r-1} = n, \omega_r = (i,j,G,\boldsymbol{\lambda})) \\
&\leq \phi_r[(\lambda_1 + \lambda_2)m^* - \varepsilon_1 N] < -\varepsilon' \phi_r N \leq -\frac{\varepsilon'}{M_0} N
\end{aligned} \quad (16)$$

as $\phi_r = \max(a_r, b_r) \geq 1/M_0$ and $\lambda_1 + \lambda_2 \leq 2M_0$ by Condition E.

Now we apply Theorem 5.1 to $Z_n = f(\omega, r_n, \xi_n)$, with $C$ chosen so that $\tau_A$, the time $\{\xi_n\}$ hits $\mathbf{B}_A$, satisfies $\tau_A < \zeta_C$ (the time $\{Z_n\}$ falls to level $C$). This is possible because, by Lemma 4.2, $f(\omega, r, \mathbf{x}) \geq \|x\|/M_0$, that is, $\|x\|$ large implies $f$ large. Hence,

$$\mathbf{E}_\omega(\tau_A \mid \xi_0 = \mathbf{x}) < \frac{f(\omega, 0, \mathbf{x}) M_0}{\varepsilon'} < \infty$$

for each $\omega \in F$ and this in turn implies that the process $\xi$ enters $\mathbf{B}_A$ infinitely often along almost all parameter sequences. Our process is not irreducible in the normal sense, as it is not possible to reach the empty state $\varnothing$ in epoch $r$ from state $(r, \mathbf{x})$, where $\mathbf{x} = (x,y)$ with $y > 0$. However, it is possible to reach $\varnothing$ in epoch $r+1$ and the probability of this is bounded away from 0 uniformly in $\omega$ (the bound does depend upon $A$ and $M_0$). Hence, by Borel–Cantelli, we have $\tau < \infty$ a.s.

(ii) Here $\sum \pi_i L_i = 0$. Again, we use $f$ to denote the remaining time for the vector field model, but this time with the natural field $\mathbf{D}^0$. We have from Lemma 4.2 that $f$ is finite a.s. and from the inequality (12), we know that for each parameter sequence $\omega \in F$ the process $(f(\omega, r_n, \xi_n))_n$ is a supermartingale (remember that we only consider instantaneous switching in this case). Now, from standard results, it follows that $\tau_A$ is a.s. finite for every $\omega \in F$. That $\tau$ is a.s. finite follows as in part (i).

(iii) We have $\sum \pi_i L_i > 0$ and we wish to show that $\tau = \infty$ with positive probability. Consider the initial value function $g$ for a vector field $\mathbf{D}^-$ uniformly below $\mathbf{D}^0$ but with $\sum_i \pi_i L_i^- > 0$. We know from Lemma 3.3(iii) and Lemma 4.2(ii) that, for almost all parameter sequences $\omega$, the $x$-coefficients $c_r$ of $g$ in epoch $r$ converge to 0 as $r \to \infty$. For such $\omega$ and for any given constants $A > 0$ and $K' > 0$, there exists $r_0$ such that [recall (15)] the set $(r, \mathbf{B}_A) \subset \mathcal{G}_{K'}(\omega)$ for all epochs $r \geq r_0$. For any fixed parameter sequence $\omega$, the process $\xi$ is certain to reach some state $(r', \mathbf{x}') \notin \mathcal{G}_{K'}(\omega)$ where $r' \geq r_0$. Let $K = g(\omega, r', \mathbf{x}')$. Starting from $(r', \mathbf{x}')$, $\rho_{K'} \leq \tau_A$ a.s. by construction but, as we now show, $\rho_{K'} < \infty$ a.s. is impossible. The conditions of Lemma 5.2 hold from which we get $E_\omega(U_\infty) \leq 1/K$, but if $\rho_{K'} < \infty$ a.s., then $U_\infty \geq 1/K' > 1/K$ a.s. As $U_\infty \geq 0$, it follows that $P_\omega(\tau_A = \rho_{K'} = \infty) \geq 1 - K'/K > 0$ for almost all $\omega$.

REMARK. Another way to describe the idea in part (iii) of the proof above is to say when $(r, \mathbf{B}_A) \subset \mathcal{G}_{K'}$ for all epochs $r$ there is some probability $\beta > 0$ (uniform over parameter sequences $\omega$) such that $\mathbf{P}_\omega(\xi_n \notin \mathcal{G}_{K'}$ for all $n) > \beta$, that is, if the vector field $\mathbf{D}^-$ has trajectories that remain outside $\mathcal{G}_{K'}$, then the random walk $\xi$ can stay outside $\mathcal{G}_{K'}$ indefinitely.

5.4. *Semimartingale estimates for Theorem* 2.2(i). We now state the key result we use for establishing finiteness of moments of hitting times to the neighborhood of the empty state. It is a slight extension of Theorem 1 of [1]. This, in turn, is the generalization to general powers $s$ of Theorem 2.2 of [7] (a similar result for integer moments only).

THEOREM 5.3. *Let $K > 1$, $(\Omega, (\mathcal{F}_n)_n, \mathbf{P})$ be a filtered probability space and $(Z_n)_{n \in \mathbb{N}}$ an $(\mathcal{F}_n)_n$-adapted positive valued process with $Z_0 = z \geq K$. Let $\sigma_K = \inf\{n \geq 1 : Z_n \leq K\}$, $\sigma'$ an arbitrary $(\mathcal{F}_n)_n$-stopping time and $\sigma = \inf\{\sigma_K, \sigma'\}$. Suppose that there exists $s_0 > 0$ such that, for every $n \in \mathbb{N}$, we have $\mathbf{E}(Z_n^{s_0}) < \infty$. If there exists $\varepsilon_2 > 0$ such that*

$$\mathbf{E}(Z_{n+1}^{s_0} - Z_n^{s_0} \mid \mathcal{F}_n) \leq -\varepsilon_2 Z_n^{s_0 - 1} \qquad \text{on } \{\sigma > n\},$$

*then, for all $s \leq s_0$, there exists $C = C(s_0, \varepsilon_2)$ such that $\mathbf{E}(\sigma^s \mid Z_0 = z) \leq C z^s$.*

REMARK. If we try to apply this result to $Z_n = f(\omega, r_n, \xi_n)$, where $f$ is the remaining life of the natural field $\mathbf{D}^0$, the condition of the theorem is not generally satisfied. This is why we have to introduce the modified field $\mathbf{D}^+$.

PROOF OF THEOREM 5.3. Straightforward extension (to deal with the introduction of $\sigma'$) of the corresponding proof of Theorem 1 of [1]. □

The following result establishes the key inequality for the application of Theorem 5.3 under a condition on $\eta^0(s)$, the Perron–Frobenius eigenvalue of $M^0(s)$ which was considered in Lemma 4.4. As in the proof of Theorem 2.1(i), it is necessary to sample the process at carefully chosen times to deal with the absence of service during switching times.

PROPOSITION 5.4. *Suppose that $\eta^0(s_0) < 1$ for some $s_0 > 0$ and let $\mathbf{D}^+$ be uniformly above $\mathbf{D}^0$ and $s_0$-neutral. Write $Z_n = f(\omega, r_n, \xi_n)$ for all $\omega \in F$, where $f$ denotes the remaining life under $\mathbf{D}^+$. Then for such $\omega$ and all $s < s_0$, there exists $\varepsilon_2 > 0$ such that at times $n$ which are not switching times and for $\|\mathbf{x}\|$ sufficiently large,*

$$\mathbf{E}_\omega(Z_{n+1}^s - Z_n^s \mid r_n = r, \xi_n = \mathbf{x}) \leq -\varepsilon_2 Z_n^{s-1}.$$

PROOF. Letting $z = Z_n = f(\omega, r, \mathbf{x})$, we have, for $\|\mathbf{x}\|$ (and, hence, $z$) large, by Lemma 4.2, $1/C\|\mathbf{x}\| \leq \phi_r/z \leq C/\|\mathbf{x}\|$ for some $C > 0$ which is uniform in $\omega$. Writing $\Delta f_n = f(\omega, r_{n+1}, \xi_{n+1}) - f(\omega, r_n, \xi_n)$, we have

$$\mathbf{E}_\omega(Z_{n+1}^s - Z_n^s \mid r_n = r, \xi_n = \mathbf{x}) = \mathbf{E}_\omega((z + \Delta f_n)^s - z^s \mid r_n = r, \xi_n = \mathbf{x})$$

and from Taylor's theorem,

$$(z + \Delta f_n)^s = z^s + s\Delta f_n z^{s-1} + s(s-1)\int_0^1 (\Delta f_n)^2 (z + t\Delta f_n)^{s-2}(1-t)\,dt.$$

As $\Delta f_n \geq -\phi_r$, when $r_n = r$, we have, for $t \in [0,1]$,

$$(z + t\Delta f_n)^{s-2} \leq z^{s-2}\left(1 - \frac{\phi_r}{z}\right)^{s-2} \quad \text{if } s \leq 2$$

and so, from Lemma 4.5,

$$\mathbf{E}(Z_{n+1}^s - Z_n^s \mid r_n = r, \xi_n = \mathbf{x}) \leq \phi_r\left(-s\varepsilon_1 + C_1 s(s-1)\frac{\phi_r}{z}\right) z^{s-1},$$

where $C_1 > 0$ is some constant which is uniform in $\omega$.

For $s > 2$, we have

$$(z + t\Delta f_n)^{s-2} \leq \begin{cases} z^{s-2}(1+t)^{s-2}, & \text{if } \Delta f_n \leq z, \\ (\Delta f_n)^{s-2}(1+t)^{s-2}, & \text{if } \Delta f_n > z. \end{cases}$$

The assumption of a uniform bound $K_s$ on moment $s$ of the service time distribution implies there is some constant $C'$ such that $\mathbf{E}_\omega(\Delta f_n^s \mid r_n = r, \xi_n = \mathbf{x}) \leq \phi_r^s C' K_s$ and, hence, $\mathbf{E}_\omega(\Delta f_n^s I\{\Delta f_n > z\} \mid r_n = r, \xi_n = \mathbf{x}) = O(\phi_r^s)$. In this case, for $\|\mathbf{x}\|$ large enough, there exists some $C_1 > 0$ (uniform in $\omega$) such that

$$\mathbf{E}(Z_{n+1}^s - Z_n^s \mid r_n = r, \xi_n = \mathbf{x}) \leq \phi_r z^{s-1}\left(-s\varepsilon_1 + C_1 s(s-1)\frac{\phi_r}{z} + O((\phi_r/z)^{s-1})\right),$$

again using Lemma 4.5.

Finally, for any given $s > 2$ and any $\varepsilon_2$ with $0 < \varepsilon_2 < s\varepsilon_1$, we can choose $\|\mathbf{x}\|$ large enough that $-s\varepsilon_1 + C_1 s(s-1)\frac{\phi_r}{z} + O((\phi_r/z)^{s-1}) < -\varepsilon_2$. The argument is similar for $s \leq 2$. This establishes the result as $\phi_r \geq 1/M_0$. □

PROOF OF THEOREM 2.2(i). As in Proposition 5.4, we consider the process $f(\omega, r_n, \xi_n)$ where $f$ is the remaining life under a vector field $D^+$ which is uniformly above $\mathbf{D}^0$ and $s$-neutral. It follows from $\mathbf{E}(f(\omega, r, \mathbf{x})^{s_0}) < \infty$ (Lemma 4.5) and the fact that $\mathbf{E}_\omega((\xi_{n1} + \xi_{n2})^s) < \infty$ for all $\omega$, $s > 0$ and all $n$ that $\mathbf{E}(f(\omega, r_n, \xi_n)^s) < \infty$ for all $s \in [0, s_0]$ and all $n$.

To deal with switching times, we sample this process at times $N_n$ where $N_{n+1} - N_n = 1$ unless $N_n$ is the end of an epoch, in which case $N_{n+1} - N_n = N$ for some constant $N$ to be determined. This makes it necessary to introduce modified versions of the hitting times for which we wish to make moment estimates.

Let $Z_n = f(\omega, r_{N_n}, \xi_{N_n})$, $\tilde{\sigma}_K = \inf\{n \geq 1 : Z_n \leq K\}$, $\tilde{\tau}_A = \inf\{n \geq 1 : \xi_{N_n} \in \mathbf{B}_A\}$ for suitably large $A$ and $\tilde{\sigma} = \min\{\tilde{\sigma}_K, \tilde{\tau}_A\}$. For $Z_0 = z > K$ for some sufficiently large constant $K$, we have, by Proposition 5.4,

$$(17) \qquad \mathbf{E}_\omega(Z_{n+1}^s - Z_n^s \,|\, r_n = r, \xi_n = \mathbf{x}) \leq -\varepsilon_2 Z_n^{s-1} \qquad \text{on } \{\tilde{\sigma} > n\}$$

as long as jump $N_n$ is not a switching event. If the server switches queues at time $N_n$, then

$$Z_{n+1} = Z_n + f(\omega, r_{N_n} + 1, \xi_{N_n+N} - \xi_{N_n}).$$

By the argument leading to the inequality (16), we have, for large enough $N$ (but $\xi_{N_n} > N$ to ensure there is no change of epoch),

$$\mathbf{E}(f(\omega, r_{N_n} + 1, \xi_{N_n+N} - \xi_{N_n}) \,|\, \mathcal{F}_{N_n}) \propto -N.$$

Using this together with the moment bound assumed above, it follows that, for given $\varepsilon_2$, there exists $N$ such that the bound (17) also holds when the jump at time $N_n$ is a switching event.

Hence, by Theorem 5.3, there exists a constant $C$ such that $\mathbf{E}_\omega(\tilde{\sigma}^s \,|\, Z_0 = z) < Cz^s$ uniformly over $\omega \in F$. We can choose $A > KM_0$ which ensures that $\tilde{\tau}_A \leq \tilde{\sigma}_K$ a.s. and, hence, for $\mathbf{x} \notin \mathbf{B}_A$,

$$\mathbf{E}(\tilde{\tau}_A^s \,|\, \xi_0 = \mathbf{x}) = \mathbf{E}(\tilde{\sigma}^s \,|\, Z_0 = z)$$
$$= \mathbf{E}(\mathbf{E}_\omega(\tilde{\sigma}^s \,|\, Z_0 = z))$$
$$\leq C\mathbf{E}(f(\omega, 1, \mathbf{x})^s) < \infty.$$

To complete the proof, let $\tau_A = \inf\{n \geq 1 : \xi_n \in \mathbf{B}_A\}$ for suitably large $A$ (so $\tau_A$ is measured in process time $n$ rather than the sampling sequence $\{N_n\}$) and recall that $\tau = \inf\{n \geq 1 : \xi_n = \varnothing\}$. Observe that $\tau_A < N\tilde{\tau}_A$ so that $\mathbf{E}(\tau_A^s \,|\, \xi_0 = \mathbf{x}) < \infty$. Finally, by a standard argument, for example, Theorem A.1 from [10], it follows that $\mathbf{E}(\tau^s \,|\, \xi_0 = \mathbf{x}) < \infty$ for all finite $\mathbf{x}$. □

5.5. *Proof of Theorem* 2.2(ii), *nonexistence of moments.* We start this proof by showing that, with a probability that can be bounded away from 0 uniformly for all parameter sequences, the random walk started outside $\mathcal{G}_K$ stays outside $\mathcal{G}_K$ at least until it hits $\mathbf{B}_A$.

PROPOSITION 5.5. *Suppose that $\eta^0(s_0) < 1$ for some $s_0 > 0$ and that $g$ is the initial value function of an s-neutral vector field $\mathbf{D}^-$. For any given parameter sequence $\omega$, suppose that the queue length process $\xi$ starts at state* $(r, \mathbf{x})$ such that $g(\omega, r, \mathbf{x}) = K$ and $\mathbf{x} \notin \mathbf{B}_A$ for $A$ large enough to satisfy the conditions of Lemma 5.2. For $1 < K' < K$,

$$\mathbf{P}_\omega(\rho_{K'} \geq \tau_A) \geq 1 - \frac{K'}{K} > 0.$$

PROOF. As $\eta^0(s_0) < 1$, Theorem 2.1(iii) implies that $\tau_A$ is a.s. finite. Consider the process $U_n = 1/g(\xi_{n \wedge \tau_A \wedge \rho_{K'}})$. By Lemma 5.2, $(U_n)$ is a positive super-martingale. Therefore, the limit $U_\infty$ exists almost surely and its law is supported by $[0, 1]$ as $K' > 1$. Moreover, by Fatou's lemma,

$$\frac{1}{K} = \mathbf{E}_\omega(U_0) \geq \mathbf{E}_\omega(U_\infty)$$
$$= \frac{1}{K'}\mathbf{P}_\omega(\rho_{K'} < \tau_A) + \mathbf{E}_\omega(U_{\tau_A} \mid \rho_{K'} \geq \tau_A)\mathbf{P}_\omega(\rho_{K'} \geq \tau_A)$$

as $\tau_A$ is a.s. finite. Hence, $\mathbf{P}_\omega(\rho_{K'} < \tau_A) \leq \frac{K'}{K} < 1$ by the initial choice of $K$ and $K'$ and $(U_n)$ being nonnegative. Thus, $\mathbf{P}_\omega(\rho_{K'} \geq \tau_A) \geq 1 - \frac{K'}{K} > 0$. □

PROOF OF THEOREM 2.2(ii). We have $\eta^0(s) > 1$ and $\eta^0(s_0) < 1$ for some $s_0 < s$. Hence, the remaining life $f$ of vector field $\mathbf{D}^-$ (uniformly below $\mathbf{D}^0$) is a.s. finite. As $f$ and the initial value function $g$ of $\mathbf{D}^-$ are both linear within epochs, $g$ can be scaled so that $g(\omega, 1, x\mathbf{e}_1) = f(\omega, 1, x\mathbf{e}_1)$ for all $x > 0$. Starting $\xi$ from state $(r, \mathbf{x})$ with $g(\omega, r, \mathbf{x}) = K > K'$ and applying Proposition 5.5, we see that $\mathbf{P}_\omega(\rho_{K'} \geq \tau_A) \geq 1 - K'/K > 0$ uniformly over all $\omega$.

We now show that on the event $\{\rho_{K'} \geq \tau_A\}$ the random walk takes such a long time to reach $\mathbf{B}_A$ that the $s_0$ moment of $\tau_A$ cannot be finite. The event $\{\rho_{K'} \geq \tau_A\}$ can only happen when (i) there is an epoch $\hat{r}$ where $(\hat{r}, \mathbf{B}_A)$ is not a subset of $\mathcal{G}_{K'}$, that is, where the $\mathbf{D}^-$ trajectory $\{g(\omega, S(t), V(t)) : t \geq 0, \ g(\omega, S(t), V(t)) = K'\}$ enters $\mathbf{B}_A$ during epoch $\hat{r}$; (ii) the random walk does not enter $\mathcal{G}_{K'}$ prior to epoch $\hat{r}$.

For fixed $\omega$ and $K'$, suppose the random walk starts from $\xi_0 = (x_0, 0)$ in epoch 1 with $g(\omega, 1, x_0\mathbf{e}_1) > K'$. Let the random times $T_1, T_2, \ldots$ denote the ends of the epochs (for the random walk, not the dynamical system) up to epoch $\hat{r} - 1$. On the event $\{\rho_{K'} \geq \tau_A\}$ we have

$$f(\omega, 1, x_0\mathbf{e}_1) \leq \sum_{r=1}^{\hat{r}-1} t_r - t_{r-1} + AM_0/m_0,$$

where, in the notation of Section 3.2, the $t_r$ are the epoch endtimes for a particle moving according to the vector field $\mathbf{D}^-$ and the bound $AM_0/m_0$ on the final epoch's contribution to $f(\omega, 1, x_0\mathbf{e}_1)$, the time for the particle trajectory to reach $\mathbf{B}_A$, is obtained from Condition E with some obvious geometry. Recalling the notation $v_1(\omega, t_{r-1})$ for the particle position at the start of epoch $r$ (for the dynamical system) and writing $X_{r-1}$ for the queue length at the start of epoch $r$ (for the random walk), we have, for each epoch $r < \hat{r}$,

$$T_r - T_{r-1} \geq X_{r-1} > v_1(\omega, t_{r-1}),$$

as only one job can complete in each service time and, by assumption, the random walk does not enter $G_{K'}$. From

$$v_1(\omega, t_{r-1}) = d_1(t_r - t_{r-1}) \geq m_0(t_r - t_{r-1}),$$

we have

$$\tau_A > T_{\hat{r}-1} = \sum_{r=1}^{\hat{r}-1} T_r - T_{r-1} > \sum_{r=1}^{\hat{r}-1} v_1(\omega, t_{r-1}) \geq m_0 \sum_{r=1}^{\hat{r}-1} t_r - t_{r-1}$$
$$\geq m_0 f(\omega, 1, x_0 \mathbf{e}_1) - AM_0$$

on all random walk trajectories such that $\tau_A \leq \rho_{K'}$.

Finally, as $\eta^-(s) > 1$, for any initial vector $\mathbf{x} \notin B_A$ for $A$ sufficiently large,

$$\mathbf{E}(\tau_A^s \mid \xi_0 = \mathbf{x}) \geq (1 - K'/K) m_0^s \mathbf{E}([f(\omega, 1, K'\mathbf{e}_1) - AM_0]^s) = \infty,$$

where the final equality follows from Lemma 3.4 applied to the vector field $\mathbf{D}^-$. $\square$

## REFERENCES

[1] ASPANDIIAROV, S., IASNOGORODSKY, R. and MENSHIKOV, M. (1996). Passage-time moments for nonnegative stochastic processes and an application to reflected random walks in a quadrant. *Ann. Probab.* **24** 932–960. MR1404537
[2] DE CALAN, C., LUCK, J. M., NIEUWENHUIZEN, TH. M. and PETRITIS, D. (1985). On the distribution of a random variable occurring in 1D disordered systems. *J. Phys. A* **18** 501–523. MR0783195
[3] FAYOLLE, G., MALYSHEV, V. A. and MENSHIKOV, M. V. (1995). *Topics in the Constructive Theory of Countable Markov Chains*. Cambridge Univ. Press. MR1331145
[4] FAYOLLE, G., IGNATYUK, I., MALYSHEV, V. and MENSHIKOV, M. (1991). Random walks in two-dimensional complexes. *Queueing Systems* **9** 269–300. MR1132178
[5] FOSS, S. and LAST, G. (1996). Stability of polling systems with exhaustive service policies and state-dependent routing. *Ann. Appl. Probab.* **6** 116–137. MR1389834
[6] LAMPERTI, J. (1960). Criteria for the recurrence or transience of stochastic processes. I. *J. Math. Anal. Appl.* **1** 314–330. MR0126872
[7] LAMPERTI, J. (1963). Criteria for stochastic processes. II. Passage-time moments. *J. Math. Anal. Appl.* **7** 127–145. MR0159361
[8] MACPHEE, I. M. and MENSHIKOV, M. (2003). Critical random walks on two-dimensional complexes with applications to polling systems. *Ann. Appl. Probab.* **13** 1399–1422. MR2023881
[9] MACPHEE, I., MENSHIKOV, M., PETRITIS, D. and POPOV, S. YU. (2007). Polling systems with random changes of regime. *Ann. Appl. Probab.* To appear.

[10] Menshikov, M. V. and Popov, S. Yu. (1995). Exact power estimates for countable Markov chains. *Markov Process. Related Fields* **1** 57–78. MR1403077
[11] Menshikov, M. and Zuyev, S. (2001). Polling systems in the critical regime. *Stochastic Process. Appl.* **92** 201–218. MR1817586
[12] Meyn, S. P. and Tweedie, R. L. (1993). *Markov Chains and Stochastic Stability.* Springer, London. MR1287609

I. MacPhee
M. Menshikov
Mathematics Department
University of Durham
United Kingdom
E-mail: i.m.macphee@durham.ac.uk
        mikhail.menshikov@durham.ac.uk

D. Petritis
Insitut de recherche mathématique
Université de Rennes 1
France
E-mail: dimitri.petritis@univ-rennes1.fr

S. Popov
Instituto de matemática e estatística
Universidade de São Paulo
Brasil
E-mail: popov@ime.usp.br